\documentclass[12pt]{article}

\usepackage[english]{babel}
\usepackage{amsmath}
\usepackage{amssymb,amsfonts,amsthm}
\usepackage{tikz,color,pgf,graphicx,url}
\usetikzlibrary{calc}
\usepackage{enumerate}     % več možnosti za sezname

\usepackage{xcolor} % for colors

\usepackage{hyperref}
\hypersetup{
	colorlinks   = true, %Colours links instead of ugly boxes
	urlcolor     = black, %Colour for external hyperlinks
	linkcolor    = black, %Colour of internal links
	citecolor   = black %Colour of citations
}

\usepackage[
top=3cm,
bottom=3cm,
inner=2.5cm,      % margini za dvostransko tiskanje
outer=2.5cm,
footskip=40pt     % pozicija številke strani
]{geometry}

\newtheorem{theorem}{Theorem}[section]
\newtheorem{corollary}[theorem]{Corollary}
\newtheorem{lemma}[theorem]{Lemma}
\newtheorem{proposition}[theorem]{Proposition}
\newtheorem{conjecture}[theorem]{Conjecture}

\newtheorem{question}[theorem]{Question}
\newtheorem{definition}[theorem]{Definition}

\newtheorem{example}[theorem]{Example}

\definecolor{lavenderindigo}{rgb}{0.58, 0.34, 0.92}

\DeclareMathOperator{\lcapt}{lcapt}

\begin{document}
	
	\title{Localization game capture time of trees and outerplanar graphs}
	
	\author{
		Vesna Ir\v si\v c Chenoweth $^{a,b,}$\thanks{Email: \texttt{vesna.irsic@fmf.uni-lj.si}}
		\and
		Matija Skrt $^{a,}$\thanks{Email: \texttt{matija.skrt@gmail.com}}
	}
	
	\maketitle
	
	\begin{center}
		$^a$ Faculty of Mathematics and Physics, University of Ljubljana, Slovenia\\
		\medskip
		
		$^b$ Institute of Mathematics, Physics and Mechanics, Ljubljana, Slovenia\\
		\medskip
	\end{center}
	
	\begin{abstract}
		The localization game is a variant of the game of Cops and Robber in which the robber is invisible and moves between adjacent vertices, but the cops can probe any $k$ vertices of the graph to obtain the distance between probed vertices and the robber. The localization number of a graph is the minimum $k$ needed for the cops to be able to locate the robber in finite time. The localization capture time is the minimum number of rounds needed for the cops to locate the robber. 
		
		The localization capture time conjecture claims that there exists a constant $C$ such that the localization number of every connected graph on $n$ vertices is at most $Cn$. While it is known that the conjecture holds for trees, in this paper we significantly improve the known upper bound for the localization capture time of trees. We also prove the conjecture for a subclass of outerplanar graphs and present a generalization of the localization game that appears useful for making further progress towards the conjecture.
	\end{abstract}
	
	
	\noindent
	{\bf AMS Subj.\ Class.\ (2020)}: 05C57, 05C12, 05C05
	
	%%%%%%%%%%%%%%%%%%%%%%%%%%%%%%%%%
	%%%%%%%%%%%%%%%%%%%%%%%%%%%%%%%%%
	\section{Introduction}
	\label{sec:intro}
	%%%%%%%%%%%%%%%%%%%%%%%%%%%%%%%%%
	%%%%%%%%%%%%%%%%%%%%%%%%%%%%%%%%%
	
	The Cops and Robber game on graphs has been widely studied; see for example \cite{book-1}. The game is played on a graph by two players, one controlling the cops and the other the robber. Positions of the cops and the robber are known to both players at all times. The players are allowed to move to an adjacent vertex or stay on the same vertex, and players alternate taking moves. The goal of the cops is to capture the robber, i.e.\ at least one cop must occupy the same vertex as the robber. The robber's goal is to evade capture forever. The minimum number of cops needed to ensure cops winning on a graph $G$ is the \emph{cop number} $c(G)$ of $G$. In 2009, Bonato et al.\ \cite{bonato+2009} introduced the capture time of a graph as the minimum number of moves needed for $c(G)$ cops to win on $G$. The problem is well-understood for cop-win graphs, but only a handful of results are known for graphs with a larger cop number \cite{kinnersley-2018, pisantechakool+2016, brandt+2020}.
	
	In 2012, Seager \cite{seager-2012} introduced a version of Cops and Robber where the robber is invisible and the cops are locating the robber using probes. This was later refined to the game called the \emph{localization game} which has been further studied for example in \cite{bonato+23-diam2, bonato-2021, bonato+20, bonato+2024, bosek+18-geometric-planar, boshoff+21, brandt+2017, brandt+2020, carraher+2012, haslegrave+2018}. This game is also played on a graph by two players, one controlling the cops and the other the robber. However, the robber is invisible to the cops during the game, while the robber knows the positions of the cops. The players alternate taking moves, with the robber selecting the starting vertex first. The robber is allowed to stay on the same vertex or move to a neighbor. The cops are allowed to move to any set of vertices (not necessarily adjacent to their previous positions), and after they all move, each cop sends out a cop probe, which gives the distance between each cop to the robber. Thus after every move of the cops, they obtain a distance vector $(d(c_1, r), \ldots, d(c_k,r))$ where $c_1, \ldots, c_k$ are positions of the cops and $r$ is the position of the robber. The goal of the cops is to determine the position of the robber after finitely many rounds. The robber's goal is to evade capture forever. Additionally, we assume that the robber is omniscient, meaning that the robber is aware of all strategies of the cops. The \emph{localization number} $\zeta(G)$ of $G$ is the minimum number of cops needed to ensure that the cops win the localization game on $G$. We only consider this game for connected graphs. The exact upper bound for the localization number is known for example for trees, outerplanar graphs \cite{bonato+20}, and some diameter 2 graphs \cite{bonato+23-diam2}, and is known to be unbounded for planar graphs \cite{bosek+18-geometric-planar}.
	
	In 2022, Behague et al.\ \cite{behague+22} introduced in full generality the analogous problem to the capture time for the Cops and Robber game for the localization game. The \emph{$k$-localization capture time} $\lcapt_k(G)$ is the minimum number of rounds needed for $k \geq \zeta(G)$ cops to capture the robber on $G$ if the cops are minimizing the number of rounds and the robber is maximizing it. If $k = \zeta(G)$, we refer to the parameter as the \emph{localization capture time} and denote it by $\lcapt(G)$. Recall that the game starts by the robber choosing an initial vertex. Then, round 1 begins. Each round consists of a move of the cops, after which they obtain a distance vector to the robber's position (we also say that the cops \emph{probe} the $k$ vertices they move to), which is followed by a move of the robber. The number of rounds needed for the cops to determine the position of the robber is thus equal to the number of times the cops obtain a distance vector to the robber, i.e.\ the number of times the cops probe the vertices. Note that as $k$ increases, $\lcapt_k(G)$ decreases, and it equals 1 if $k$ is the metric dimension of $G$.
	
	As mentioned already in \cite{behague+22}, determining the localization capture time of graphs is extremely challenging. They determine several bounds for the $k$-localization capture time if $k > \zeta(G)$ and provide an upper bound for the localization capture time of incidence graphs of projective planes. They also initiate the study of a general upper bound for the localization capture time in terms of the number of vertices of the graph. A family of graphs is \emph{well-localizable} if there exists a constant $D$ such that for every graph $G$ in the family, $\lcapt(G) \leq D |V(G)|$. They conjecture the following.
	
	\begin{conjecture}[{Localization Capture Time Conjecture (LCTC), \cite{behague+22}}]
		\label{conj:LCTC}
		The family of all connected graphs is well-localizable.
	\end{conjecture}
	
	Moreover, they wonder if an even stronger result holds.
	
	\begin{question}[\cite{behague+22}]
		\label{que:LCTC}
		If $G$ is a connected graph on $n$ vertices, then $\lcapt(G) \leq n$.
	\end{question}
	
	It is known that trees, interval graphs and complete multipartite graphs are well-localizable \cite{behague+22}. In particular, for graphs on $n$ vertices, they show that $\lcapt(T) \leq 5n$ if $T$ is a tree with $\zeta(T) = 1$, $\lcapt(T) \leq n-1$ if $T$ is a tree with $\zeta(T) = 2$, $\lcapt(G) \leq n$ if $G$ is an interval graph, and  $\lcapt(G) \leq n-1$ if $G$ is a complete multipartite graph.
	
	In this paper, we improve the best known bounds for the localization capture time for trees, in particular proving that Question \ref{que:LCTC} holds true for all trees. We also show that outerplanar graphs $G$ with $\zeta(G)=2$ are well-localizable, and moreover, that Question \ref{que:LCTC} holds true for them as well. As proving Conjecture \ref{conj:LCTC} seems extremely difficult in general, especially for families of graphs for which the exact localization number is not known (which is most families of graphs), we devote the rest of the paper to presenting a generalization of the localization game related to colorings of graphs. We believe that this more general approach could be used to prove general upper bounds for the localization capture time.

	\section{Trees}
	\label{sec:trees}
	
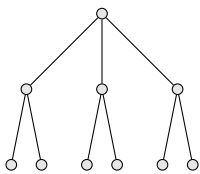
\begin{figure}
		\centering
		\begin{tikzpicture}
			[scale=1,
			vert/.style={circle, draw, fill=black!10, inner sep=0pt, minimum width=4pt}, 
			double/.style={circle, draw, fill=black, inner sep=0pt, minimum width=4pt}, 
			central/.style={circle, draw, fill=black, inner sep=0pt, minimum width=4pt},
			]
			
			\node[vert] (1) at (0,0) {};
			\node[vert] (2) at (-1,-1) {};
			\node[vert] (3) at (0,-1) {};
			\node[vert] (4) at (1,-1) {};
			\node[vert] (5) at (-1.2,-2) {};
			\node[vert] (6) at (-0.8,-2) {};
			\node[vert] (7) at (-0.2,-2) {};
			\node[vert] (8) at (0.2,-2) {};
			\node[vert] (9) at (0.8,-2) {};
			\node[vert] (10) at (1.2,-2) {};
			
			\draw (5) -- (2) -- (1) -- (4) -- (10);
			\draw (7) -- (3) -- (8);
			\draw (6) -- (2);
			\draw (9) -- (4);
			\draw (3) -- (1);
		\end{tikzpicture}
		\caption{Tree $T_{3,3}$.}
		\label{fig:T33}
	\end{figure}

	We start this section with a few definitions. As in \cite{seager-2014}, let $T_{3,3}$ be the tree in Figure \ref{fig:T33}. Let $\ell(T)$ be the number of leaves of a tree $T$. Recall that the set of vertices of a rooted tree can be partitioned into several levels, level $i$ being the set of vertices of distance $i$ from the root. Additionally, if $v$ is a vertex in a rooted tree, then we call the set containing a child of $v$ and all its descendants a \emph{branch}. Note that $v$ can have multiple branches.
	
	It is known that one cop wins the localization game on a tree $T$ if and only if $T$ does not contain a copy of $T_{3,3}$ (\cite[Theorem 8]{seager-2014}). As mentioned in \cite{bosek+18-centroidal} and explicitly proven in \cite{bonato-book}, it is known that if $T$ is a tree, then $\zeta(T) \leq 2$. It follows from \cite{seager-2014} and \cite{behague+22} that the best known upper bound for the localization capture time for a tree $T$ on $n$ vertices is 
	$$\lcapt(T) \leq \begin{cases}
		5n; & \text{if } \zeta(T) = 1,\\
		n-1; & \text{if } \zeta(T) = 2.
	\end{cases}$$
	In this section we improve on these two upper bounds. First observe that $\lcapt(K_1) = \lcapt(K_2) = 1$ and that $\lcapt(K_{1, m}) = m-1 = |V(K_{1,m})|-2$ for $m \geq 2$. As $\ell(T) \leq n-2$ for every tree on $n \geq 3$ vertices that is not a star, and $\ell(T) \leq n-4$ for every tree on $n$ vertices that contains a copy of $T_{3,3}$, we get the following result from Theorems \ref{thm:tree-1} and \ref{thm:tree-2} which are proved below.
	
	\begin{corollary}
		\label{cor:trees}
		If $T$ is a tree on $n\geq3$ vertices, then 
		$$\lcapt(T) \leq \begin{cases}
			n-2; & \text{if } \zeta(T) = 1,\\
			\left \lfloor \frac{n}{2} \right \rfloor -3; & \text{if } \zeta(T) = 2.
		\end{cases}$$
	\end{corollary}
	
	As already mentioned, the bound $n-2$ cannot be improved for trees $T$ with $\zeta(T) = 1$. After proving Theorem \ref{thm:tree-2}, we also provide an example showing that the bound is tight for trees $T$ with $\zeta(T) = 2$. In the proofs, the candidate set denotes the set of all possible locations of the robber on a graph after a move by the cops or the robber is completed.
	
	\begin{theorem}
		\label{thm:tree-1}
		If $T$ is a tree on at least two vertices that does not contain a copy of $T_{3,3}$, then $\lcapt(T) \leq \ell(T)$.
	\end{theorem}
	\begin{proof}
		Let $T$ be a tree which does not contain a copy of $T_{3,3}$, and root $T$ at any leaf $v_0$. We
		provide a strategy for the cop to construct a sequence of pairwise distinct vertices $v_1, v_2, \ldots, v_a$ of $T$, 
		lying on levels $k_1 < k_2 < \cdots < k_a$ respectively, such that for each $1 \leq i \leq a$, the candidate set will be contained in the set of descendants of $v_i$,
		immediately after the cop's probe. 
		In particular, after each application of one of the strategies described below, the cop will be able to either capture the robber, or identify the lowest 
		positioned vertex $v_i$, on level $k_i > k_{i-1}$,  
		such that the candidate set will be contained in the set of its descendants throughout the remainder of the game. 
		Since the tree is finite, the cop will indeed capture the robber in a finite number of moves.
		Moreover, we prove that during the described strategy, the cop (in most cases) maintains the following property: the number of probes needed to eliminate a part of the tree from all the possible future candidate sets with respect to the outlined strategy is at most the number of leaves in that part of the tree. There are a few exceptions to this general rule which are dealt with on the spot.
		
		Choose $v_0$ as the first probe. If the robber is not captured, then the candidate set is clearly contained in the set of descendants of $v_0$. Hence, there exists a unique vertex $v$, positioned lowest out of all the vertices for which the candidate set is contained in the set of their descendants. We set $v_1 = v$. Suppose that after some move, the cop has already 
		constructed $v_1,\ldots, v_t$.  Since the robber must be located on a descendant of $v_t$, we may now consider the levels of these vertices with respect to $v_t$,
		that is, as if the tree was rooted at $v_t$, in order to simplify the notation. Recall that the vertex set consisting of a child of $v_t$ and all of its descendants is called a branch.
		
		\begin{description}
			\item[Case 1.] $v_t$ has at most one child of degree 3 or more.
			
			Suppose that after a probe by the cop, the robber is known to be on level $k>0$, which we temporarily fix. Recall that the levels are considered with respect to $v_t$.
			
			First, the cop consecutively probes all of the leaves adjacent to $v_t$. If he ever gets the answer $0$ or $1$, the robber is captured. 
			Thus, after completing this procedure, the cop knows, that the robber must be located on one of the other branches. 
			If there are no leaves incident to $v_t$, the cop simply omits this step. Note that to eliminate each such leaf one probe is needed, which is exactly equal to the number of leaves in the eliminated part.
			
			Suppose now that there is at least one child of $v_t$ of degree 2. Again, if such children are not present, the cop simply omits this step.
			Otherwise, he consecutively probes the children of such vertices.
			Let $u$ be one such vertex that the cop probes, and let $d$ be the answer he gets after this probe. 
			Immediately before the probe, the robber cannot be on any other level than $k-1, k$ or $k+1$. 
			We consider the following four cases. In all of them, if the cop ever realizes that the robber is on the same branch as him and on level 3 or lower, we may simply set $v_{t+1}$ to be
			$u$ or, if possible, some other vertex on the branch on a lower level.
			
			\begin{description}
				\item[Case 1.a.] $d=0$\\
				The robber is captured.
				
				\item[Case 1.b.] $d=1$\\
				We consider two possibilities. If there exists a branch that the cop has not probed on yet, the cop probes $v_t$ (here, we consider all different branches from $v_t$).
				If he gets an answer of 0, 1 or 2, the robber is captured. 
				Otherwise, the robber is located on the branch on which the cop probed on the previous move, and on level 3 or lower. Notice that since the cop has not probed on at least one of the other branches, he has indeed not used too many moves, since the unprobed branch contains at least one leaf.
				
				Suppose now, that the cop has already probed all of the other branches (if there is no other branch to probe, the robber is known to be located on level $k$). 
				If the robber was on level 1 or 2, he would have been captured on the previous move. Hence, without any further probes, the cop knows that the robber is located on the last branch, on level three or lower.
				
				\item[Case 1.c.]  $d=2$\\
				If $k=1$, the cop knows that the robber is on $v_t$, thus capturing the robber. 
				If $k \in \{3, 4, 5\}$, the robber is on the same branch as the cop and is located on level 3 or lower. Note that no other value of $k > 0$ is compatible with $d=2$ at this stage of the game. 
				
				\item[Case 1.d.]  $d \geq 3$\\
				If the robber is located on the same branch as the cop, the possible values of $d$ are $k-3,k-2,$ and $k-1$. If he is located on a different branch, the possible values of $d$ are $k+1, k+2,$ and $k+3$. Hence, if $d<k$, the robber is located on the same branch as the cop, on level 3 or lower, and if $d>k$, the robber is located on a different branch. In the second case, the cop chooses another child of $v_t$ of degree 2, and repeats the process.
			\end{description}
			
			Finally, if the robber is not found on a leaf incident to $v_t$ and also not on a branch which begins with a vertex of degree 2,
			the cop knows, after the last probe, that the robber is either located on the child of $v_t$ on the unprobed branch, thus capturing the robber, or is contained in the set of 
			all of the descendants of that child. Hence, we can set $v_{t+1}$ to be
			the child of $v_t$ on the unprobed branch or, if possible, some other vertex on the branch on a lower level. 
			
			\item[Case 2.] $v_t$ has two children of degree 3 or more. 
			
			Note that  $v_t$ cannot have three or more such children as it does not contain $T_{3,3}$.
			
			We call the two branches that contain the children of $v_t$ of degree 3 or more the \emph{main} branches.
			The goal of the cop's strategy will be to reduce the case where there are 
			two main branches to the case where there is at most one. If the cop successfully identifies the branch on which the robber is located, or 
			makes sure that he is on a branch, but not on at least one of the main branches, in a small enough amount of moves,
			he may then use the strategy described above, to extend the sequence $(v_i)$. Notice that this reduction is perfectly legitimate,
			since the strategy we have described above makes sure that the robber cannot reach $v_t$ without getting captured.
			
			Let $d_i$ be the distance the cop receives as the answer to his $i$-th probe (after the selection of $v_t$).
			We call the main branch with more leaves the \emph{rich} branch, and the other 
			main branch the \emph{poor} branch. If they both have the
			same amount of leaves, we make the choice between them arbitrarily. Notice that since we have assumed 
			that there are two main branches, this means that $v_t$ is his parent's only child, since
			$T$  does not contain a copy of $T_{3,3}$ and we rooted the tree at a leaf.
			
			Note that the strategy of the cop is to eliminate the rich or the poor branch (or both), thus reducing Case 2 to Case 1, and then to continue using the strategy described in Case 1.
			
			\item[Case 2.1.]The rich branch contains at least three leaves.
			
			We now distinguish
			several cases, depending on the level on which the robber is located. Recall that after the last probe, the robber is known to be on level $k$ in the subtree of $T$ consisting of $v_t$ and all its descendants.
			
			\begin{description}
				\item[Case 2.1.a.] $k=1$\\
				The cop first probes the child of $v_t$ on the poor branch. If $d_1=0$, the robber is captured. If $d_1 \geq 2$, the cop knows that the robber is on some other branch.
				Therefore, assume $d_1=1$. The cop now probes the child of $v_t$ on the rich branch. If $d_2=0$ or $d_2=1$, the robber is captured. 
				If $d_2 \geq 3$, the robber must be on the poor branch. 
				Suppose now, that $d_2=2$. 
				If $v_t$'s parent is not $v_0$, then the cop probes the parent of the parent of $v_t$. If $d_3 \in \{0,1,2\}$, the robber is captured. Otherwise $d_3 \geq 3$ and the robber is located on some branch, but
				not on the rich branch. Similarly, if $v_t$'s parent is $v_0$, the cop probes $v_0$ and obtains the same conclusions, with $d_3$ being reduced by 1. In any case, at most three probings were needed to eliminate the rich branch which has at least three leaves by assumption.
				
				\item[Case 2.1.b.] $k=2$\\
				This time, the cop first probes the child of $v_t$ on the rich branch. If $d_1=0$, the robber is captured. If $d_1 =1$, the robber is located on the rich branch. If
				$d_1 \geq 3$, the robber is on some other branch. Suppose therefore that $d_1 = 2$. This implies, that the robber is either located on the rich branch on level 3, or on some other branch on level 1. 
				The cop now probes the child of $v_t$ on the poor branch. If $d_2=0$, the robber is captured. 
				If $d_2 \geq 2$, the robber is located on a branch other than the poor branch. Therefore, assume $d_2=1$. Analogously to the previous case, the cop now probes on the child of $v_t$ on the rich branch, and then 
				on the parent of the parent of $v_t$, thus either capturing the robber, or making sure that he is on a branch other than the rich one. Since the cop might have needed four probes to eliminate the rich branch (which could only have three leaves), we must analyze this case further. 
				
				If the rich branch contains at least four leaves, the cop has not used too many moves. However, the same is also true if the rich branch contains three leaves.
				Suppose that this is the case, and let $p$ be the number of leaves on the poor branch. 
				Clearly, $p \in \{2, 3\}$. If the candidate set is ever contained in the poor branch, the cop 
				simply probes $p-1$ of its leaves consecutively, capturing the robber in at most $p-1$ moves,
				recovering from the move he had lost before. In particular, if $p=2$, probing on one of the leaves uniquely determines the position of the robber. If $p=3$, there are two branches stemming from the child of $v_t$ on the poor branch, where one of them must have at most one leaf. Probing on this solitary leaf either captures the robber or reduces to the case $p=2$ using only one move. Otherwise, the cop locates the robber on a different branch, 
				thus eliminating the poor branch, and giving himself at least two extra moves to spare.
				
				\item[Case 2.1.c.] $k \geq 3$\\
				The cop again first probes the child of $v_t$ on the rich branch. If $d_1 \geq k+1$, the cop knows that the robber is located on some other branch.
				If $d_1 \leq k-1$, the robber is located on the rich branch.
				Therefore assume $d_1=k$, that is, the robber is either located on level $k+1$ on the rich branch, or on level $k-1$ on some other branch. 
				The cop now probes the child of $v_t$ on the poor branch. 
				If $d_2 \geq k$, the robber is located on a branch, but not on the poor branch. If $d_2 \leq k-1$, the robber is located on a branch, but not on the rich branch. Here, at most two probings were needed to eliminate one of the main branches, which both have at least two leaves by assumption.
				
			\end{description}
			
			\item[Case 2.2.] Both of the main branches contain exactly two leaves.\\
			The cop first probes one of the leaves
			on the rich branch. If the robber is located on the path from this leaf to the child of $v_t$ on this branch, he is captured. 
			Otherwise, the cop now probes the other leaf of the branch. If the robber is not captured, he is either located on some other branch, or perhaps on the parent of $v_t$. 
			The cop now simply probes the parent of the parent of $v_t$, or $v_0$ if that is the parent of $v_t$. The robber is either captured,
			or located on one of the branches other than the rich branch. The cop has not used too many moves for the following reason.
			If the candidate set is ever contained in the poor branch, the cop needs at most one move to capture the robber (as in Case 2.1.b) while the poor branch has two leaves. Otherwise,
			the cop has eliminated the poor branch too, gaining two extra moves for free.
		\end{description}
		
		Notice that in each case, the cop has eliminated some branches as the robber's possible whereabouts, using at most as many moves as there are 
		leaves on the discarded branches.
		Hence, using the described strategy, the cop extends his sequence of vertices, until the robber is 
		captured, needing at most $\ell(T)$ moves to do so.
	\end{proof}

    We do not know if the upper bound in Theorem \ref{thm:tree-1} is tight. However, many trees $T$, including paths and stars, have capture time equal to $\ell(T)-1$.
	
	\begin{theorem}
		\label{thm:tree-2}
		If $T$ is a tree that contains a copy of $T_{3,3}$, then $\lcapt(T) \leq \left\lfloor\frac{\ell(T)}{2}\right\rfloor - 1$.
	\end{theorem}
	\begin{proof}
		Let $T$ be a tree as stated. By \cite[Theorem 8]{seager-2014} we have $\zeta (T) = 2$. We use a simplified version of the cop's strategy provided in the proof of Theorem \ref{thm:tree-1}, making sure that we remove at least two new leaves from the set of vertices on which the robber might be located in the rest of the game with each move, and removing at least four in the first move. Hence, our strategy will use at most $\lfloor\frac{\ell(T)}{2}\rfloor - 1$ moves, as desired.
		
		We begin by rooting $T$ at the vertex $v_0$ that has at least three children $x, y, z$, which each have at least two children. This is possible since $T$ contains a copy of $T_{3,3}$. For our first move we probe on $x$ and $y$. If the distances obtained are distinct, the robber is located on one of the probed branches, and we set $v_1$ to be on the respective branch. If not, and there are only three branches at $v_0$, set $v_1$ to be the appropriate vertex on the unprobed branch, i.e.\ $v_1 = z$ or a lower vertex if possible. Otherwise, there are at least two more branches to be checked, and we set $v_1 = v_0$, ignoring the branches that we have already discarded. Notice that we have discarded at least two branches, each having at least two leaves, from the candidate set. We continue as follows.
		
		Suppose that after some move, the cop has already constructed $v_0,\ldots, v_t$.  Since the robber must be located on a descendant of $v_t$, we may now consider the levels of vertices to be with respect to $v_t$, that is, as if the tree was rooted at $v_t$. Clearly, by the choice of $v_t$, we must have at least two branches stemming from $v_t$. We now simply divide the branches into pairs, possibly with one remaining unpaired, and in each move we probe the neighbors of $v_t$ on the respective pair of branches. If the robber is located on any of the two, we identify it, and set $v_{t+1}$ to be the child of $v_t$ on the respective branch, or a lower positioned vertex, if that is possible, as in the proof of Theorem \ref{thm:tree-1}. Assuming that there is still at least one branch left unprobed, we have eliminated at least two branches, each containing at least one leaf. Otherwise, we move on to the next pair. Note that the robber is captured if he ever moves to $v_t$ since that is the unique vertex for which the distances to both of the cops are equal to one. 
		
		Suppose now, that we have discarded all of the branches, except for the last two or three (depending on whether the number of branches is even or odd). Then the robber must be located on these final two or three branches. If there are three, we simply probe the next pair of branches and thus locate the robbers' branch. Note that we have eliminated at least two branches in this move. Hence, assume that there are only two. If they both contain only one leaf, we capture the robber in the next move, by probing on both of those leaves, thus uniquely determining the robber's position. If they both contain more than one leaf, we simply probe both of them on the neighbors of $v_t$ and move on as usual -- discarding at least one of the branches. The problem arises when exactly one of the branches has exactly one leaf. Then we do the following. Let $u$ be the first descendant of $v_t$ on the remaining branch with at least two leaves that has at least two children. If $u$ has exactly two children, we probe both of them in the next move. With that we are successfully able to determine whether the robber is contained in one of the two branches stemming from $u$ (if the obtained distances are distinct), and set $v_{t+1}$ accordingly, or capture the robber otherwise. If, however, $u$ has more than two children, we simply probe on the first vertex of the branch stemming from $v_t$ with exactly one leaf and the neighbor of $u$ on one of the branches stemming from it. If the robber was located on $v_t$, on the branch with one leaf, or on a vertex between $v_t$ and $u$, he can be captured. Otherwise, we know that the robber is not located on at least one of the branches stemming from $u$ nor on the branch stemming from $v_t$ that has only one leaf, and may continue the described procedure with the appropriately set $v_{t+1}$.
	\end{proof}
	
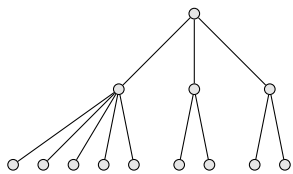
\begin{figure}
		\centering
		\begin{tikzpicture}
			[scale=1,
			vert/.style={circle, draw, fill=black!10, inner sep=0pt, minimum width=4pt}, 
			double/.style={circle, draw, fill=black, inner sep=0pt, minimum width=4pt}, 
			central/.style={circle, draw, fill=black, inner sep=0pt, minimum width=4pt},
			]
			
			\node[vert] (1) at (0,0) {};
			\node[vert] (2) at (-1,-1) {};
			\node[vert] (3) at (0,-1) {};
			\node[vert] (4) at (1,-1) {};
			\node[vert] (5) at (-1.2,-2) {};
			\node[vert] (6) at (-0.8,-2) {};
			\node[vert] (7) at (-0.2,-2) {};
			\node[vert] (8) at (0.2,-2) {};
			\node[vert] (9) at (0.8,-2) {};
			\node[vert] (10) at (1.2,-2) {};
			\node[vert] (11) at (-1.6,-2) {};
			\node[vert] (12) at (-2,-2) {};
			\node[vert] (13) at (-2.4,-2) {};
			
			\draw (5) -- (2) -- (1) -- (4) -- (10);
			\draw (7) -- (3) -- (8);
			\draw (6) -- (2);
			\draw (9) -- (4);
			\draw (3) -- (1);
			\draw (2) -- (11);
			\draw (2) -- (12);
			\draw (2) -- (13);
		\end{tikzpicture}
		\caption{Tree $G_5$.}
		\label{fig:T33-extended}
	\end{figure}

	Note that this upper bound also cannot be improved in general. Let $G_m$ be a tree obtained from $T_{3,3}$ and adding $m-2$ leaves to one of the children of the root. For example, see Figure \ref{fig:T33-extended}. Then we need at least one move to eliminate the part of the tree that does not contain the star $K_{1,m}$, and $\left\lfloor\frac{m}{2}\right\rfloor$ moves for the star part, giving $\left\lfloor\frac{m}{2}\right\rfloor + 1 = \left\lfloor\frac{m+4}{2}\right\rfloor - 1 = \left\lfloor\frac{\ell(G_m)}{2}\right\rfloor - 1 = \left\lfloor\frac{|V(G_m)|}{2}\right\rfloor - 3$.
	
	\section{Outerplanar graphs}
	\label{sec:outerplanar}
	
	Bonato and Kinnersley showed in \cite{bonato+20} that $\zeta(G) \leq 2$ for every outerplanar graph $G$. In this section we extend their results by considering the localization capture time of outerplanar graphs in the game with two cops. Our results prove that outerplanar graphs with localization number 2 are well-localizable. It remains to consider outerplanar graphs with localization number 1 to fully answer Question \ref{que:LCTC} for all outerplanar graphs.
	
	Recall that a \emph{block} of a graph is a maximal connected subgraph that has no cut-vertex. The next result relies strongly on the winning strategy given in \cite{bonato+20}. For clarity, we only point out the differences and additions made to their results. We recall a few basic definitions.
	
	The two cops playing on an outerplanar graph maintain a \emph{cop territory} during the game which is a set of vertices with specific properties. Two distinct vertices of the cop territory are the \emph{endpoints}, sometimes referred to as $v_L$ and $v_R$. The cop territory is a connected subgraph of the outerplanar graph, after every probe the cops know that the robber is not on a vertex from the cop territory, no vertex (except possibly the endpoints) of the cop territory is adjacent to a vertex outside the cop territory, and both endpoints belong to the same block.
	
	\begin{proposition}
		\label{prop:outerplanar}
		If G is an outerplanar graph that is an edge-disjoint union of blocks $B_1, \ldots, B_m$, then $\lcapt_2(G) \leq \sum_{i=1}^m (|B_i| - 1)$.
	\end{proposition}
	\begin{proof}
		We consider the winning strategy for two cops on an outerplanar graph $G$ described in the proof of \cite[Theorem 3.1]{bonato+20}. In that strategy, at least one vertex is added to the cop territory in each move, except in the Case 4(b). If $w_2 = v_{r+1}$, we would have to use two moves to add just one vertex to the cop territory. However, we may easily fix this by skipping the probe on $v_s$ and $w_2$ and probe instead $v_L$ and $w_2$ directly, thus saving a move.
		
		Notice also that we add two vertices to the cop territory in the first move. Furthermore, each time we move into a new block, we must already have at least one of its vertices contained in the cop territory. Hence, we use at most $|B_i| - 1$ moves in every block, and the result follows.
	\end{proof}
	
	Note that Proposition \ref{prop:outerplanar} immediately implies that if $G$ is an outerplanar graph on $n$ vertices with $\zeta(G) = 2$, then $\lcapt(G) \leq n-1$ and thus $G$ satisfies Question \ref{que:LCTC}. However, we can obtain an even better upper bound for 2-connected outerplanar graphs.
	
	\begin{proposition}
		\label{prop:outerplanar-2-con}
		If G is a 2-connected outerplanar graph with $c$ chords, then $\lcapt_2(G) \leq c + 1$.
	\end{proposition}
	
	\begin{proof}    
		If $c=0$, then $G$ is a cycle and one move clearly suffices. Suppose $c \geq 1$. Order the vertices of the outer face of $G$ in the clockwise direction as $v_1, \ldots, v_n$, such that $\deg(v_n) \geq 3$. For $1 \leq i<j\leq n$, let $v_i-v_j := G[v_i, v_{i+1}, \ldots, v_j]$. Let now $C = v_nv_1\ldots v_av_n$ be an induced cycle in $G$. We probe on $v_{\lfloor \frac{a+1}{2} \rfloor - 1}$ and $v_{\lfloor \frac{a+1}{2} \rfloor}$. The robber is either captured or the candidate set is contained in $v_a-v_n$, in which case we set $v_L = v_n$, and $v_R = v_a$. Note, that we also know the distances from $v_L$ and $v_R$ to the robber, due to the probe.
		
		Similarly as in the proof of \cite[Theorem 3.1]{bonato+20} we consider two cases. In each one, we choose vertices $v_L$ and $v_R$ with $L > R$, where each of them lies on a chord contained in $v_R-v_L$, such that (after the procedure described in each case is completed) either the robber is captured, or we know the distances $d_R$ and $d_L$ from the robber to the vertices $v_R$ and $v_L$ respectively, and the candidate set is contained in $v_R-v_L$.
		
		\begin{description}
			\item[Case 1.] $v_R$ and $v_L$ do not have common neighbors in $v_R-v_L$.\\
			Suppose without loss of generality that $d_R \geq d_L$ and let $w$ be the (clockwise) furthest neighbor of $v_R$ in $v_R-v_L$. Then we can take the new endpoints $v_L$ and the unique vertex $v$ in $w-v_L$ that is closest to $w$ and lies on a chord contained in $v-v_L$. Note that if no such $v$ existed, the robber would have been caught on the previous move. The new distance $d_R$ can then be obtained by subtracting the distance between $v_R$ and $v$ from $d_R$.
			
			\item[Case 2.] $v_R$ and $v_L$ have a common neighbor in $v_R-v_L$.\\
			If the robber is closer to one endpoint, we do the same as in Case 1. Suppose therefore that $d_L = d_R$ and let $v_s$ be the common neighbor of $v_R$ and $v_L$ in $v_R-v_L$. We use a modified version of the argument from Case 4(b) in the proof of \cite[Theorem 3.1]{bonato+20}, discarding at least one chord with each move. Let $v_R = w_1, w_2, \ldots, w_k$ be the neighbors of $v_s$ in $G$ that are counterclockwise from $v_s$. For $i \in \{1, \ldots, k-1\}$ let sector $i$ be $w_i-w_{i+1}$. The cops aim to determine which sector (if any) the robber occupies. We again consider two cases.
			
			\begin{itemize}
				\item[\textbf{a)}] $v_R-w_2v_sv_R$ is an induced cycle and $v_L$ and $v_R$ are not neighbors.\\
				Let $v$ be a vertex on the arc $v_R-w_2$, such that $d(v_R, v) = d(w_2, v) + \{0, 1\}$. The cops now probe on $v_L$ and $v$, obtaining distances $d_L'$ and $d_v$. If the robber is deduced to be contained in the interior of $v_R-w_2$, another round finishes the game. If not, we may safely take $v_L$ and $w_2$ as the new $v_L$ and $v_R$, knowing both $d_L$ and $d_R$ (from $d_L'$ and $d_v$).
				
				\item[\textbf{b)}] $v_R-w_2v_sv_R$ is not an induced cycle or $v_L$ and $v_R$ are neighbors.\\
				In this case, the cops probe $v_s$ and $w_2$; let $d_s$ and $d_w$ denote the distances observed. If $d_s \geq d_w$, then the robber cannot presently reside in $v_s-v_L$, since every shortest path from $w_2$ to a vertex in $G$ must pass through either $v_L$ or $v_s$, and $v_s$ is closer to both of these than $w_2$ is. In this case, as before, the cops may take $v_s$ and $v_R$ as the endpoints. Thus we may suppose that $d_s < d_w$; since $v_s$ and $w_2$ are adjacent, we must have $d_w = d_s +1$.
				We claim that the robber cannot occupy sector 1. 
				
				Suppose to the contrary that the robber does occupy some vertex in sector 1, and note that it cannot be $w_2$ (since the cops have just probed $w_2$). Since $d_w = d_s +1$, some shortest path from $w_2$ to the robber passes through $v_s$ and, since the robber is in sector 1, through $v_R$ as well. We thus have $d_L \in d_w + 2 + \{-1, 0, 1\}$ and $d_R \in d_w - 2 + \{-1, 0, 1\}$, and so $d_w+1 \leq d_L = d_R \leq d_w-1$, a contradiction.
				
				We now set $v_R = w_2$, leaving $v_L$ unchanged and probe $v_R$ and $v_L$. Note that we can afford to spare two moves on discarding sector 1, since we made sure that in doing so, we discarded at least 2 chords.             
			\end{itemize}
			
			We now repeat the process for $i = 3, \ldots, k-1$, where the role of $v_R$ is played by $w_i$. If we find ourselves in the case $i=k-1$ with $d_R = d_L$ and establish that the robber is not contained in sector $k-1$, we do not set $v_R = w_k$, but $v_R = v_s$, as we trivially can do. This finishes the case, $v_L$ and $v_s$ are now the proper endpoints of the cop territory, and $d_L$ and $d_R$ are known, and we may therefore move on to the next case.
			
			It is clear that the candidate set is contained in $v_R-v_L$ throughout, since the robber being located on $v_L$, $v_R$ or $v_s$ results in his immediate capture.
			
		\end{description}
		
		In each case, $v_R-v_L$ (in which the candidate set is contained) contains at least as many chords less than before the endpoints were changed, as is the number of rounds necessary to complete a specific case. Taking the preparatory round 1 into account, the robber is captured in at most $c+1$ moves.
	\end{proof}
	
	Let $G$ be an outerplanar graph on $n$ vertices with $c$ chords and with $\zeta(G) = 2$. As $G$ is outerplanar, $|E(G)| \leq 2n-3$ and $|E(G)|$ equals the sum of the number of edges on the outer cycle ($n$) and the number of chords ($c$). Thus $c \leq n-3$. Hence, Proposition \ref{prop:outerplanar-2-con} immediately implies that $\lcapt(G) \leq n-2$, so $G$ satisfies Question \ref{que:LCTC}, and the bound is always better than the one given by Proposition \ref{prop:outerplanar}.
	
	\section{Generalization}
	\label{sec:generalization}
	In this section we introduce a generalization of the game. After giving the definition we also formally prove that it is indeed a generalization of the localization game (see Theorem \ref{thm:equivalence}) and explain why this more general approach is promising to make progress towards Conjecture \ref{conj:LCTC}.

    For $S \subseteq V(G)$ and a coloring $C$ of $V(G)$, we denote by $C(S)$ the restriction of $C$ to $S$ which is also a coloring. Recall that a color class of a coloring is a subset of vertices that are assigned the same color.
	
	We are given a graph $G$ and colorings $C_1, \ldots, C_k$ of $V(G)$. We build the \emph{game structure} of $G$ according to $C_1, \ldots, C_k$ starting with row 1 and adding new rows consecutively. 
	\begin{itemize}
		\item Row 1 consists of singletons $\{v\}$ for every $v \in V(G)$.
		\item Row $i$ consists of all subsets $S$ of $V(G)$ which do not appear in rows $1, \ldots, i-1$ for which there exists $j \in [k]$ such that the color classes of $C_j(N[S])$ appear in rows $1, \ldots, i-1$.
	\end{itemize}
	
	A set $S\subseteq V(G)$ is \emph{monochromatic} in a coloring $C_i$ if all of its vertices are of the same color in $C_i$.
	
	To every set $S\subseteq V(G)$ that lies in the game structure in row $i \geq 2$, we adjoin a set of colorings $C_j$, such that the color classes of $C_j(N[S])$ appear in rows $1, \ldots, i-1$.

    If $S$ is in row $i+1$ and $T$ is in row $i$, then $S$ is a \emph{parent} of $T$ and $T$ is a \emph{child} of $S$ if there exists a coloring adjoined to the set $S$, in which $T$ is monochromatic.
	
	From the game structure we obtain the \emph{reduced game structure} as follows. In the highest row we keep only the sets which have all subsets lying in lower rows. Then, moving on to the highest but one row, we only keep the sets whose subsets all lie in lower rows and that have a parent in the (newly obtained) highest row. We remove the other sets from this row.
	
	Next, we apply the same procedure on one row lower (with respect to the newly obtained highest but one row) to discard some of the sets. Continuing until we reach row 1 we obtain the desired reduced structure. By definition, every set in row $i+1$ has a child in row $i$. Taking a minimal such set (with respect to inclusion) shows that the reduced structure has the same amount of rows as the original one, i.e.\ no rows have become empty during the procedure.

	\begin{example}
        \label{ex:4.1.}
		Let $T$ be the tree in Figure \ref{fig:example-tree} and let colorings of $V(T)$ be obtained from comparing distances to a given vertex. Then the obtained color classes are
		\begin{align*}
			C_1 &: \{ 1,23, 45 \},\\
			C_2 &: \{ 1,2,3, 45 \},\\
			C_3 &: \{ 145,2,3 \},\\
			C_4 &: \{ 15,2,3,4 \},\\
			C_5 &: \{ 14,2,3,5 \}.
			%C_6 &: \{ 145,2,3,6 \}.
		\end{align*}
		
		\begin{figure}
			\centering
			\begin{tikzpicture}
				[scale=1,
				vert/.style={circle, draw, fill=black!10, inner sep=0pt, minimum width=4pt}, 
				double/.style={circle, draw, fill=black, inner sep=0pt, minimum width=4pt}, 
				central/.style={circle, draw, fill=black, inner sep=0pt, minimum width=4pt},
				]
				
				\node[vert, label=above:$1$] (1) at (0,0) {};
				\node[vert, label=left:$2$] (2) at (-1,-1) {};
				\node[vert, label=right:$3$] (3) at (1,-1) {};
				\node[vert, label=below:$4$] (4) at (0.5,-2) {};
				\node[vert, label=below:$5$] (5) at (1.5,-2) {};
				%\node[vert, label=below:$6$] (6) at (1.5,-2) {};
				
				\draw (2) -- (1) -- (3) -- (4);
				\draw (5) -- (3); %-- (6);
			\end{tikzpicture}
			\caption{Tree $T$.}
			\label{fig:example-tree}
		\end{figure}
		
		Then the obtained game structure is:
		\begin{align*}
			\text{Row 3: }& 13,23,34,35,123,134,135,145,234,235,245,345,1234,1235,1245,1345,2345,12345;\\
			\text{Row 2: }& 12,14,15,24,25,45,124,125;\\
			\text{Row 1: }& 1,2,3,4,5.
		\end{align*}
		
		The reduced game structure then becomes:
		\begin{align*}
			\text{Row 3: }& 13,23,34,35,145,245;\\
			\text{Row 2: }& 14,15,45;\\
			\text{Row 1: }& 1,2,3,4,5.
		\end{align*}
	\end{example}
	
	Given a graph $G$ on $n$ vertices and a positive integer $k$, the set of \emph{distance colorings} is obtained by taking every $S \subseteq V(G)$, $|S| = k$, and defining the coloring $C_S \colon V(G) \to [n]^k$ as $C_S(v) = (d(v,s_1), \ldots, d(v, s_\ell))$ where $S = \{s_1, \ldots, s_\ell\}$. Note that this set of colorings is exactly the one we used in Example \ref{ex:4.1.} with $k = 1$.
	
	\begin{theorem}
		\label{thm:equivalence}
		Let $G$ be a graph, let $k \geq \zeta(G)$ be a positive integer, and let $C_1, \ldots, C_m$ be the distance colorings of $G$. Then the number of rows obtained in this structure is $lcapt_k(G)+1$.
	\end{theorem}
	
	\begin{proof}
		The following correspondences between the game structure and the localization game hold: 
		\begin{itemize}
			\item each set $S$ in the game structure corresponds to the possible whereabouts of the robber after a probe by the cops,
			\item the set $N[S]$ corresponds to the candidate set after some round is completed,
			\item the adjoined colorings of $S$ correspond to the probes that the cops may use in order to bring the robber one move closer to capture.
		\end{itemize}
		
		Keeping this in mind, notice also that due to $k \geq \zeta(G)$, every nonempty subset of $V(G)$ appears in the game structure. Row 2 contains all subsets $S$ of $V(G)$ for which there exists a probe on some other vertex set $T$, such that if the candidate set is ever equal to the (closed) neighborhood of $S$, the cops may probe on $T$, corresponding to the coloring $C_T$ in which all vertices of the neighborhood of $S$ are colored distinctly (that is, the cops differentiate between these vertices and can thus determine robber's location). In other words, row 2 contains precisely those vertex sets such that if their neighborhoods are a candidate set, the cops capture the robber in the next round. Analogously, the row $i+1$ contains those vertex sets $S$, whose neighborhoods, should they become a candidate set, can be probed with the probe corresponding to a respective adjoined coloring of $S$, after which the cops know for sure that the robber is located on one of the sets in the lower rows. Hence, inductively, the row $i>1$ contains precisely those vertex sets, such that if their neighborhoods should become candidate sets, the robber can be captured in $i-1$ moves. Thus clearly the height of the game structure is at most $\lcapt_k(G)+1$.
		
		We similarly provide a strategy for the robber. In the beginning, the candidate set is simply $V(G)$ which is located in the highest row, due to the above considerations. Suppose the candidate set in round $j$ is the neighborhood of a set $S$ in row $i$. Whichever probe the cops choose in round $j+1$, the robber may always make sure that the candidate set is the neighborhood of some set lying on row $i-1$ or higher -- if the cops probe on $T$, he chooses as his ``position'' (he need not choose a single vertex since he knows the moves of the cops in advance) the highest lying color class of $C_T(N[S])$. Hence $\lcapt_k(G) + 1$ must be greater than or equal to the height of the game structure (the robber is already captured on row 1).
	\end{proof}

    If we wish to also include colorings corresponding to the cops probing on vertices which are not necessarily distinct, we can relax the condition $|S| = k$ to $|S| \leq k$. But as long as $k \leq |G|$, the cops need never choose such a probe in their optimal strategy, hence such colorings are irrelevant for our discussion.

    \begin{example}
        \label{ex:4.2.}
		Let $T$ be the tree $T_{3, 3}$ in Figure \ref{fig:example-two-cops} and let colorings of $V(T)$ be distance colorings obtained from two cops.
		
		\begin{figure}
			\centering
			\begin{tikzpicture}
				[scale=0.7,
				vert/.style={circle, draw, fill=black!10, inner sep=0pt, minimum width=4pt}, 
				double/.style={circle, draw, fill=black, inner sep=0pt, minimum width=4pt}, 
				central/.style={circle, draw, fill=black, inner sep=0pt, minimum width=4pt},
				]
                
				\node[vert, label=above:$0$, red] (0) at (0,0) {};
				\node[vert, label=left:$1$, orange] (1) at (-2,-2) {};
				\node[vert, label=left:$2$, yellow] (2) at (0,-2) {};
				\node[vert, label=left:$3$, green] (3) at (2,-2) {};
				\node[vert, label=below:$4$, blue] (4) at (-2.5,-4) {};
				\node[vert, label=below:$5$, blue] (5) at (-1.5,-4) {};
				\node[vert, label=below:$6$, lavenderindigo] (6) at (-0.5,-4) {};
                \node[vert, label=below:$7$, lavenderindigo] (7) at (0.5,-4) {};
                \node[vert, label=below:$8$, violet] (8) at (1.5,-4) {};
                \node[vert, label=below:$9$, violet] (9) at (2.5,-4) {};
				
				\draw (0) -- (1) -- (4);
				\draw (1) -- (5);
				\draw (0) -- (2) -- (6);
				\draw (2) -- (7);
				\draw (0) -- (3) -- (8);
				\draw (3) -- (9);
			\end{tikzpicture}
			\caption{Tree $T_{3, 3}$ with the distance coloring on vertices $1$ and $2$.}
			\label{fig:example-two-cops}
		\end{figure}
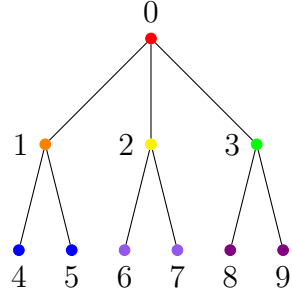
		The game structure is:
        \begin{align*}
			\text{Row 3: }& \text{Any set of vertices whose closed neighborhood contains}\\
            &\text{at least five of the vertices $4$, $5$, $6$, $7$, $8$, $9$};\\
			\text{Row 2: }& \text{Any set of vertices that is not in row $1$ or $3$};\\
			\text{Row 1: }& 0, 1, 2, 3, 4, 5, 6, 7, 8, 9.
		\end{align*}

        To see how the game structure simulates the game, let the robber choose the set of all vertices of the tree as it's starting vertex set (which is of course contained in the topmost row). This represents the starting position of the game (the candidate set is the entire graph), when the robber chooses an initial vertex, say vertex $4$. Then there exists a coloring of $T_{3,3}$ with color classes lying in rows $1$ and $2$. This is the optimal coloring (representing a probe) that the cops can choose in the game. In our case, one such coloring is the distance coloring on vertices $1$ and $2$, as shown in Figure \ref{fig:example-two-cops}. The robber can now choose any of the color classes of this coloring, and assuming optimal play, he chooses one lying in row $2$, for example $45$. This corresponds to the robber staying on vertex $4$. In the next round, the candidate set is $145$, and the cops capture the robber by choosing the distance coloring on vertices $4$ and $5$, which, of course, corresponds to a probe on vertices $4$ and $5$ in the actual game. In this way, every possible game is modeled in the game structure.
        
		The reduced game structure is then:
		\begin{align*}
			\text{Row 3: }& 123, 128, 129, 136, 137, 234, 235, 1679, 1678, 1689, 1789, 2458, 2459, 2489, 2589,\\
            &3456, 3457, 3467, 3567, 45678, 45679, 45689, 45789, 46789, 56789;\\
			\text{Row 2: }& 04, 05, 06, 07, 08, 09, 45, 67, 89, 12, 13, 23;\\
			\text{Row 1: }& 0,1,2,3,4,5,6,7,8,9.
		\end{align*}

        This reduces the number of sets in the structure, while retaining its height and the relevant parent-child relationships.
    
    \end{example}
    
	\begin{definition}
		A graph $G$ is \emph{solvable} with respect to colorings $C_1, \ldots, C_k$ of $V(G)$ if the the respective game structure contains all subsets of $V(G)$.
	\end{definition}
	
	So, a graph is solvable, if and only if it is possible to win the respective game on it, as outlined in Theorem \ref{thm:equivalence}.
	
	Observe that if $G$ is solvable for $k$ colorings, then it is also solvable for $k+1$ colorings. 
	
	Notice that this generalization gets rid of the need for calculating the distances between vertices since the colorings are given to us arbitrarily. It also generalizes the game to any subset of $V(G)$ whereas the original game only focuses on those subsets which follow from the entire set $V(G)$ according to the rules. Furthermore, in this version, the given graph and set of colorings need not be solvable, that is, the question of capture time has been fully absolved from the question of solvability. We may find out the capture time (the height of the structure) even on unsolvable combinations of graph and colorings (that is, even if we cannot cover all of the subsets with our structure). However, the number of colorings alone does not limit the height of the structure.
	
	In light of the Conjecture \ref{conj:LCTC} we pose the following.
	
	\begin{question}
		\label{que:generalization}
		Given a graph $G$ on $n$ vertices and colorings $C_1, \ldots, C_m$ of $V(G)$, does there exist a function $f$ depending only on $n$ such that the number of rows in the respective game structure is at most $f(n)$?
	\end{question}
	
	Note that a positive answer to Question \ref{que:generalization} would also provide an upper bound for $\lcapt(G)$. In particular, proving that $f(n) = C n$ for some constant $C$ suffices would imply the correctness of Conjecture \ref{conj:LCTC}.
	
	For example, if we are given only two colorings, we are able to answer Question \ref{que:generalization} positively.
	
	\begin{lemma}
		For a graph $G$ on $n$ vertices and given two colorings $C_1$ and $C_2$, the obtained game structure has at most $\binom{n}{2} + 1$ rows.
	\end{lemma}
	\begin{proof}  
		We prove that every set in the reduced structure that lies in row 3 or higher must contain exactly two vertices.
		
		By way of contradiction, suppose there exists a set $S$ in row $i \geq 3$ with $|S| \geq 3$. Since it has a parent, it must be monochromatic in one of the parent's adjoined colorings, say $C_1$. Therefore, $S$ must have $C_2$ as an adjoined coloring. Let $T$ be a child of $S$ in the reduced game structure. 
		
		We now show that $S$ is contained in $N[T]$.
		Suppose not, and let $v \in S$ be a vertex such that $N[v] \cap T = \emptyset$. Then $N[S-v]$ also contains $T$, and therefore  $S-v$ lies in the same row as $S$, a contradiction.
		
		But now there is no coloring which could be adjoined to $T$, since $S$ and $T$ are both contained in $N[T]$, while $T$ is monochromatic in $C_2$, and $S$ is monochromatic in $C_1$.
		
		Hence, each row 3 or higher contains at least one set of cardinality 2. Since at least one set of cardinality two must also be present in row 2 (row 2 would otherwise be empty), this implies the desired upper bound.
	\end{proof}

	\section*{Acknowledgments}
	
	V.I.C.\ acknowledges the financial support from the Slovenian Research and Innovation Agency (Z1-50003, P1-0297, N1-0218, N1-0285, N1-0355) and the European Union (ERC, KARST, 101071836).
	
%	\bibliographystyle{abbrv}
%	\bibliography{literatura}

\end{document}